
\magnification=\magstep1
\input amstex
\documentstyle{amsppt}
\catcode `\@=11 \def\logo@{}
\def\Monograph{\monograph@true \let\headmark\rightheadtext
  \let\varindent@\indent \def\headfont@{\bf}\def\proclaimheadfont@{\bf}
  \def\remarkheadfont@{\it}}
\def\refsfont@{\tenpoint} \catcode`\@=12
           \NoBlackBoxes
\pageheight{7.0in} \hcorrection{0.3in}
\vcorrection{0in} \baselineskip=17pt \parindent=20pt \parskip=12pt
  \define\pd#1#2{\frac{\partial #1}{\partial #2}}
\newsymbol\zR 2052\newsymbol\zZ 205A\newsymbol\zQ 2051\newsymbol\zN
204E\newsymbol \zO 204F
\font\sevenex=cmex7
\scriptfont3=\sevenex \scriptscriptfont3=\sevenex

\topmatter
\title Classification of all Poisson-Lie structures on an Infinite-Dimensional
Jet Group\endtitle
\author {Boris A. Kupershmidt} \endauthor
\affil {Department of Mathematics\\
The University of Tennessee Space Institute\\
Tullahoma, TN 37388\\
USA\\
e-mail: bkupersh\@sparc2000.utsi.edu
}
\endaffil
\endtopmatter
\topmatter
\author { Ognyan S. Stoyanov}\endauthor
\affil
{ Department of Mathematics\\
Rutgers University\\
New Brunswick, NJ 08903\\
USA\\
e-mail: stoyanov\@math.rutgers.edu
}
\endaffil
\date June 8, 1995\enddate
\endtopmatter
\NoRunningHeads

{

{\bf Abstract.} A local classification of all Poisson-Lie structures on an
infinite-dimensional group $G_{\infty}$ of formal power series is given.
All Lie bialgebra structures on the Lie algebra ${\Cal G}_{\infty}$ of
$G_{\infty}$ are also classified.
\newpage

\par Let $G_{\infty}$ be the group of formal power series in one variable
$\bigl\{x (u)=\sum_{i=1}^{\infty}x_iu^i\mid x_1\ne 0\bigr\}$, with
a group multiplication $G_{\infty}\times G_{\infty}\to G_{\infty}$ being the
substitution:
$$\bigl(xy \bigr)(u):=x \bigl(y (u)\bigr),\hskip 15pt \text{or} \hskip 15pt u
\mapsto x (u), \tag 1$$
and with an identity $e$ the identity
map $u\mapsto u$.
The group $G_{\infty}$ is the group of formal diffeomrphisms of
$\zR^1$ which leave the origin fixed. It is a
projective limit $G_{\infty}={\underset {\underset n \to\longleftarrow}
\to\lim}\,
G_n$, where $G_n$ are the
finite-dimensional Lie groups of $n$-jets of the line at the origin.
The multiplication in $G_n$ is again defined by the substitution (1):
$(\Cal X_n\Cal Y_n)(u):=\Cal X_n(\Cal Y_n(u))\mod u^{n+1}$, where
$\Cal X_n(u)$ and $\Cal Y_n(u)$ are polynomials in $u$ of degree $n$.
We define the space of smooth functions $C^{\infty}(G_{\infty})$ to be the
inductive limit $C^{\infty}(G_{\infty})={\underset {\underset n
\to\longrightarrow} \to\lim}\,C^{\infty}(G_n)$
of the spaces of smooth functions on the finite-dimensional groups $G_n$.
\par Following [1] we consider a multiplicative Poisson (Poisson-Lie)
structure on $G_{\infty}$ to be the bilinear skew-symmetric map $\{\
,\ \}\:C^{\infty}(G_{\infty})\times C^{\infty}(G_{\infty})\to
C^{\infty}(G_{\infty})$ defined by
$$\{f,g\}=\omega_{ij}\pd {f}{x_i}\pd {g}{x_j},\tag 2$$
for any $f,g\in C^{\infty}(G_{\infty})$, such
that the multiplication map $G_{\infty}\times G_{\infty}\to
G_{\infty}$ is a Poisson map. Here
$\omega_{ij}\in C^{\infty}(G_{\infty})$, for any $i,j\in\Bbb N$, and a
summation is assumed over
repeated indices.
The Poisson structure on $G_{\infty}\times G_{\infty}$ is taken to be
the product Poisson structure.
Note that the sum in (2) is
finite since by definition $f$ and $g$ are functions of finite number
of variables. Then the Jacobi identity for $\{\ ,\ \}$ implies that
$\omega_{ij}$'s satisfy
$$\omega_{ij}\pd {\omega_{kl}}{x_i}+\omega_{ik}\pd
{\omega_{lj}}{x_i}+\omega_{il}\pd {\omega_{jk}}{x_i}=0,\tag 3a$$
for any $j,k,l\in \Bbb N$. The multiplicativity of the Poisson
brackets (2) ($\{\ ,\ \}$ being a 1-cocycle) means that $\omega_{ij}$'s must
satisfy the following
infinite system of functional equations
$$\omega _{ij}(xy)=\omega _{kl}(x)\pd {z_i}{x_k} \pd {z_j}{x_l} +
\omega_{kl}{(y)}\pd {z_i}{y_k} \pd {z_j}{y_l},\qquad i,j\in\Bbb N, \tag 3b$$
where $z=xy$. Note again that the sums in the right hand side of (3b)
are finite. This is immediately seen from the explicit formulae
$$z _k =\sum_{i=1}^k
x_i\sum_{(\sum_{\alpha=1}^{i}{j_{\alpha}})=k}{y_{j_1}\ldots
y_{j_i}},\qquad k\ge 1,$$
for
the coordinates of $z$. From (3b) also follows that $\omega_{ij}(e)=0$.
\par  Do such structures exist on $G_{\infty}$?
 It is by no means obvious that such structures exist. For example:\newline
(i) Let us consider the 3-dimensional factor group
$G_3=G_{\infty}\ \text(mod)\ u^n$, for $n\ge 4$. Then there exists a
Poisson-Lie structure on $G_3$ described by
$$\align \{x_1,x_2\}&=x_1x_2\\
\{x_1,x_3\}&=4x_2^2-2x_1x_3\\
\{x_2,x_3\}&=6\frac {x_2^3}{x_1}-5x_2x_3,\endalign$$
where $x_1,x_2,x_3$ are the coordinate functions on the group $G_3$. However,
this Poisson-Lie structure can not be extended to a
Poisson-Lie structure on $G_{\infty}$.\newline
(ii) We conjecture that there are no non-trivial Poisson-Lie structures on the
group of diffeomorphisms of $S^1$.\newline
 Also a second question arises: If such structures exist, could they
be classified? The answer to the first question is given by the following
theorem.
\proclaim{Theorem 1} For every natural number $d\in \Bbb N$, and every
sequence $M_d=(\mu_n)_{n=1}^{\infty}$, such that $\mu_n=0$ for
$1\le n\le d$ and $\mu_{d+1}\ne 0$,
one has the following infinite-parameter family of Poisson-Lie structures on
$G_{\infty}$,
$$\multline\omega_{ij}(x)=\sum_{p=1}^{i}\sum_{q=1}^{j}px_pqx_q\lambda_{i-p+1,j-q+1}-\\
-\sum_{p=1}^{i}\sum_{q=1}^{j}\lambda_{pq}\sum_{\sum_{k=1}^{p}r_k=i}x_{r_1}\ldots x_{r_p}\sum_{\sum_{l=1}^{q}s_l=j}x_{s_1}\ldots x_{s_q},\endmultline\tag 4$$
where
$$\lambda_{mn}=\frac{1}{\mu_{d+1}}\biggl[\mu_{m}\lambda_{d+1,n}-\mu_{n}\lambda_{d+1,m}\biggr]\hskip
15pt \forall\ {m,n\ge 1}.\tag 5$$
Here, $\lambda_{d+1,n}$ are given by rational functions
$\lambda_{d+1,n}=\lambda_{d+1,n}(\mu_{d+1},\ldots ,\mu_{d+n})$ for ${n\ge 1}$,
which are computed by the following recursive formula
$$\lambda_{d+1,n}=-\frac
{1}{(d-n+1)\mu_{d+1}}\biggl[d\mu_{d+1}\mu_{n+d}-\sum_{s=1}^{n-1}(n+d-2s+1)\mu_{n+d-s+1}\lambda_{s,d+1}\biggr], \tag 6$$ where
$\lambda_{1,d+1}=\mu_{d+1}$, and there exists
the following single relation between the $\mu_{n}$'s (with $n\ge d+1$)
$$\mu_{2d+1}=-\frac
{1}{d\mu_{d+1}}\sum_{s=2}^{d}2(d+1-s)\mu_{2d+2-s}\lambda_{s,d+1},\tag 7$$
which are otherwise subject to no other restrictions. (We implicitly assume
that $\lambda_{mn}=0$ whenever $m<1$ or $n<1$.)
\endproclaim
The classification is given by Corollary 2 below.
\remark{Remark} The relation (7) follows from (6) when $n=d+1$, and
this is the only value of $n$ for which the expression
$(d-n+1)\lambda_{d+1,n}\mu_{d+1}$ equals zero.\endremark
\par The proof of Theorem 1 is rather technical [5,6]. We confine
ourselves to give here the main ideas and tools used. Let $\Cal
V=\{u,v,w,\ldots\}$ be a countable set and let
$C^{\infty}(G_{\infty})[[\Cal V]]$ be the ring of formal power
series generated by $\Cal V$ over $C^{\infty}(G_{\infty})$.
Let $X_i$, $i\in
\Bbb N$, be the coordinate functions on $G_{\infty}$. That is,
$X_i(x)=x_i$ for $x\in G_{\infty}$. Introduce
the formal series $X(u):=\sum_{i,j=1}^{\infty}X_iu^i$. Then
$x(u)=X(u)(x)=\sum_{i,j=1}^{\infty}x_iu^i$, and $\omega_{ij}=\{X_i,X_j\}$.
Define the formal series
$\Omega(u,v;X):=\sum_{i,j=1}^{\infty}\omega_{ij}u^iv^j$. Thus
$\Omega(u,v;X)$ is a generating series for the brackets $\omega_{ij}$.
After evaluation at $x\in G_{\infty}$ one has
$\Omega(u,v;x)=\sum_{i,j=1}^{\infty}\omega_{ij}(x)u^iv^j$. The
multiplicativity (3b) of the Poisson brackets on $G_{\infty}$ is equivalent
to $\Omega(u,v;x)$ satisfying the following functional equation
$$\Omega(u,v;xy)=\Omega(y(u),y(v);x)+\Omega(u,v;y)x'(y(u))x'(y(v)).\tag
8$$
Here $x'$ denotes the derivative  of $x$ with respect to its argument.
The general solution of (8) is given by
$$\Omega(u,v;x)=\varphi(u,v)x'(u)x'(v)-\varphi(x(u),x(v)),\tag 9$$
where $\varphi(u,v)$ is a formal series in $u,v$ subject to the
conditions:\newline
(i) $\varphi(u,v)$ is divisible by $u$ and $v$,\newline
(ii) $\varphi(u,v)=-\varphi(v,u)$.\newline
The map $\{\ ,\ \}\:C^{\infty}(G_{\infty})\times
C^{\infty}(G_{\infty})\to C^{\infty}(G_{\infty})$ induces a map $\{\
,\ \}\:C^{\infty}(G_{\infty})[[\Cal V]]\times
C^{\infty}(G_{\infty})[[\Cal V]]\to
C^{\infty}(G_{\infty})[[\Cal V]]$. In particular one has
$$\{
{X}(u),{X}(v)\}=\sum_{i,j=1}^{\infty}\{X_i,X_j\}u^iv^j={\Omega}(u,v;{X}).$$
Then the Jacobi identities (3a) are equivalent to the
single equation
$$\{ {X}(w),\{ {X}(u),{X}(v)\}\}+\{ {X}(u),\{ {X}(v),{X}(w)\}\}+\{ {X}(v),\{
{X}(w),{X}(u)\}\}=0,$$
which, after a short calculation using the explicit formula (9), implies that
$\varphi(u,v)$
must satisfy the following functional partial differential equation
$$\varphi(u,v)\bigl[
{\partial}_u\varphi(w,u)+{\partial}_v\varphi(w,v)\bigr]+c.p.=0.\tag 10$$
Thus the content of Theorem 1 is a description of all solutions of
(10) satisfying (i) and (ii).
 The relation between (9) and (4) is $\varphi
(u,v)=\sum_{i,j=1}^{\infty}\lambda_{ij}u^iv^j$. The solution (4) completely
describes the space of solutions of (10). An equivalent description can be
given as follows.
\proclaim{Theorem 1a} For each $d\in \Bbb N$, and any (formal series) $f_d(u)
,g_d(u)$ such that $f_d^{\prime}(u) g_d(u) -f_d(u) g_d^{\prime}(u) =-d\mu_{d+1}
f_d(u)$, where
$\mu_{d+1}\ne 0$ is an arbitrary parameter, and $f_d$ has a zero
of order $d+1$ at $u=0$, there is a solution of {\rm (10)} given
by
$$\varphi_d(u,v) =\frac {1}{\mu_{d+1}}\biggl[ f_d(u) g_d(v)
-f_d(v) g_d(u)\biggr].$$
The set of all solutions of {\rm (10)} is described
in this way.
\endproclaim
\proclaim{Corollary 1} A subclass of the above family is the following
countable family of Poisson-Lie structures. For each $d\in \Bbb N$,
choosing $M_d=(0,1,0,0,\ldots)$, one has
$$\multline\omega_{ij}(x)=(i-d)jx_jx_{i-d}-i(j-d)x_ix_{j-d}+\\
+x_i\sum_{\sum_{k=1}^{d+1}s_k=j}^{}x_{s_1}\ldots
x_{s_{d+1}}-x_j\sum_{\sum_{k=1}^{d+1}s_k=i}^{}x_{s_1}\ldots
x_{s_{d+1}},\endmultline\tag 11$$
for every $i,j\ge 1$. (We adopt the convention that $x_i=0$ whenever
$i<1$).\endproclaim
This family corresponds to the set of solutions
$\varphi(u,v)=uv(u^d-v^d)$, $d\in\Bbb N$, of (10).
\par Let $i_{\infty}:G_{\infty}\to G_{\infty}$ be the inversion map defined by
$i_{\infty} (x )={x}^{-1}$ for every $x\in G_{\infty}$.
\proclaim{Theorem 2} The map $i_{\infty}:G_{\infty}\to G_{\infty}$ is an
anti-Poisson map.\endproclaim
In other words one has
$\bigl\{f,g\bigr\}(i_{\infty}(x))=-\bigl\{f,g\bigr\}(x)$, for every
$f,g\in C^{\infty}(G_{\infty})$ and $x\in G_{\infty}$. The proof uses only
the explicit form (9) of the brackets on $G_{\infty}$.
Let $\overline{X}( u)$ be the inverse of $X(u)$. Then one has
$\overline{X}(X(u)) =u$, and  $X(\overline{X}(u))=u$,
as well as ${\overline{X}}{}^{\prime}(X(u)){X}^{\prime}(u)=1$,
and ${X}^{\prime}(\overline{\Cal
X}(u)){\overline{X}}{}^{\prime}(u)=1.$
On the other hand
$$\aligned 0&=\{ u,X( v)\}\\
&=\{ \overline{X}(X( u)) ,X(v )\}\\
&=\{ \overline{X}( w) ,X( v)\} {|}_{w=X( u)} +{\overline{X} }{}^{\prime}( w)
{|}_{w=X( u) }\{X( u) ,X( v)\}. \endaligned$$
Therefore,
$$\{X( v) ,\overline{X}( w)\}{|}_{w=X( u)}=\overline{X}{}^{\prime}( w)
{|}_{w=X( u)}\{X( u) ,X( v)\}. \tag 12$$
Also, one has the chain of identities
$$\aligned 0&=\{ v,\overline{X}( w)\}{|}_{w=X( u)}\\
&=\{\overline{X}(X( v)) ,\overline{X}( w)\}{|}_{w=X( u)}\\
&=\{\overline{X}( s) ,\overline{X}( w)\}{|}_{s=X( v) ,w=X(
u)}+\overline{X}{}^{\prime}( s) {|}_{s=X( v)}\{X( v) ,\overline{X}( w)\}
{|}_{w=X( u)}.\endaligned$$
Using (9) and (12), one rewrites the last identity as
$$\aligned 0&=\varphi(X(v)
,X(u))\overline{X}{}^{\prime}(X(v))\overline{X}{}^{\prime}(X(u))-\varphi(v,u)\\
&\qquad
+\overline{X}{}^{\prime}(X(v))\overline{X}{}^{\prime}(X(u))\left[\varphi(u,v){X}^{\prime}(u){X}^{\prime}(v) -\varphi(X(u) ,X(v))\right]\\
&=\{\overline{X}(s) ,\overline{X}(w)\} {|}_{s=X(v) ,w=X(u) }
+\varphi(u,v)-\overline{X}{}^{\prime}(w)\overline{X}{}^{\prime}(s)\varphi(w,s)
.\endaligned$$
Thus,
$$\{\overline{X}(w) ,\overline{X}(s)\}
=-\left[\overline{X}{}^{\prime}(w)\overline{X}{}^{\prime}(s)\varphi(w,s)
-\varphi(\overline{X}(w) ,\overline{X}(s))\right],$$
and this concludes the proof.
\remark{Remark} In the beginning of the theory of Poisson-Lie groups,
the property of the inversion map $i:G\to G$ to be anti-Poisson was
considered as an axiom [1,4]. However, for finite-dimensional groups
this property can be deduced from the other axioms [5,6]. For
infinite-dimensional groups such deduction is not likely.
 \endremark
\par To show that formula (4) provides {\it all} Poisson-Lie structures on
$G_{\infty}$, we now turn to the Lie algebra ${\Cal G}_{\infty}$ of
$G_{\infty}$. Let $\{e_n\}_{n\ge 0}$ be a basis of ${\Cal G}_{\infty}$, and let
$\alpha$ be a 1-cochain $\alpha \:{\Cal G}_{\infty}\to{\Cal
G}_{\infty}\hat\wedge {\Cal G}_{\infty}$, which
we write in the above basis as
$\alpha(e_n)=\sum_{i,j=0}^{\infty}\alpha^{ij}_{n}e_i\wedge e_j$, where
$\alpha$ takes values in the completed tensor product $\Cal
G_{\infty}\widehat{\otimes}\Cal
G_{\infty}=\bigoplus_{n=1}^{\infty}\Bigl(\oplus_{i+j=n}\Cal G_i\otimes\Cal
G_j\Bigr)$, where each $\Cal G_i$ is a one-dimensional subspace of $\Cal
G_{\infty}$ spanned by $e_i$. The Lie algebra structure on ${\Cal G}_{\infty}$
is given by
$$\bigl[e_n,e_m\bigr]=(n-m)e_{n+m}\hskip 20pt \forall\  {n,m\ge 0}.$$
Then the map $\alpha$ equips $\Cal G_{\infty}$ with a Lie bialgebra
structure [1] iff
$$\eqalign{(i)\qquad & \tau\circ\alpha=-\alpha\cr
(ii) \qquad & \alpha\bigl([a,b]\bigr)=a.\alpha(b)-b.\alpha(a),\ \ \ a,b\in\Cal
G_{\infty},\cr
(iii)\qquad & [1\otimes 1\otimes 1+(\tau\otimes
1)(1\otimes\tau)+(1\otimes\tau)(\tau\otimes 1)](1\otimes\alpha)\circ\alpha=0
,}$$
where $\tau$ is the transposition map $\tau\:\Cal
G_{\infty}\widehat{\otimes}\Cal
G_{\infty}\to\Cal G_{\infty}\widehat{\otimes}\Cal G_{\infty}$ defined by
$\tau(a\otimes b)=b\otimes a$, for any $a,b\in\Cal G_{\infty}$, and the dot
stands for the
action of $\Cal G_{\infty}$ on $\Cal G_{\infty}\widehat{\wedge}\Cal
G_{\infty}$ induced by the adjoint action of $\Cal G_{\infty}$ on itself.
In the case when $\alpha$ is a
$1$-coboundary one has $\alpha(a)=a.r$, where $r\in \Cal
G_{\infty}\hat{\wedge}\Cal G_{\infty}$
is a 0-cochain referred to as the classical
$r$-matrix [1,3]. In the latter case, (iii) above is equivalent to
$a.<~r,r>=0$, for any $a\in\Cal G_{\infty}$. Here
$<r,r>:=\bigl[r^{12},r^{13}\bigr]+\bigl[r^{12},r^{23}\bigr]+\bigl[r^{13},r^{23}\bigr]$, and
$\bigl[r^{12},r^{13}\bigr]:=
\sum_{i,j,k,l=0}^{\infty}r^{ij}r^{kl}[e_i,e_k]\wedge e_j\wedge e_l$, etc.,
where $r=\sum_{i,j=0}^{\infty}r^{ij}e_i\wedge e_j$.
\proclaim{Theorem 3} The first cohomology group $H^1({\Cal
G}_{\infty},{\Cal G}_{\infty}\hat\wedge{\Cal G}_{\infty})=~0$. That is, all
{\rm 1}-cocycles $\alpha\:{\Cal
G}_{\infty}\to{\Cal G}_{\infty}\hat\wedge{\Cal G}_{\infty}$ are
coboundaries.\endproclaim
\flushpar This can be proven by analyzing the infinite system of
linear equations
$$\alpha\bigl([e_n,e_m]\bigr)=e_n.\alpha(e_m)-e_m.\alpha(e_n),\qquad n,m\ge
0,$$
using its symmetries, and inductive arguments.
Then an analysis of the system of equations $e_n.<~r,r>=0$, $n\ge 0$,
shows that it is equivalent to $<r,r>=0$ (CYBE, the classical
Yang-Baxter equation). This turned out to be a specific property of
the algebra ${\Cal G}_{\infty}$ [5,6].
Thus, all Lie bialgebra structures on ${\Cal G}_{\infty}$ are given by
solutions of the classical Yang-Baxter equation. Moreover the
following theorem holds.
\proclaim{Theorem 4} There is a one-to-one correspondence between the
coboundary Lie bialgebra structures on ${\Cal G}_{\infty}$ and the Poisson-Lie
structures {\rm (4)} on $G_{\infty}$, the correspondence being given by
$r^{ij}=\lambda_{i+1,j+1}$, for every $i,j\ge 0$. Thus, for each $d\in\Bbb N$
we have the following infinite-parameter family of Lie bialgebra structures on
${\Cal G}_{\infty}$ given by
$\alpha(e_n)=\sum_{i,j=0}^{\infty}\alpha^{ij}_{n}e_i\wedge e_j$, where
$$\align\alpha^{ij}_{n}&=(2n-i)r^{i-n,j}+(2n-j)r^{i,j-n}\\
&=(2n-i)\lambda_{i-n+1,j+1}+(2n-j)\lambda_{i+1,j-n+1},\qquad \forall\,
{n,i,j\ge 0},\endalign$$
and $\lambda_{nm}$ are subject to the same conditions as described in Theorem
{\rm 1}.\endproclaim
\proclaim{Corollary 2} Thus, Theorem {\rm 1} describes all Poisson-Lie
structures on the group $G_{\infty}$.\endproclaim
\par The proof of Theorem 4 consists of showing that each
Lie bialgebra structure on $\Cal G_{\infty}$ can be integrated to a
unique Poisson-Lie structure on the group $G_{\infty}$. To show this
one has to show that the following infinite system of linear partial
differential equations
$$\multline\sum_{i=j}^n(i+1-j)x_{i+1-j}\pd
{\omega_{mn}}{x_i}=\omega_{m+1-j,n}(x)(m+1-j)+\omega_{m,n+1-j}(x)(n+1-j)+\\
+\sum_{k=1}^m\sum_{l=1}^n\alpha_j^{kl}(m+1-k)(n+1-l)x_{m+1-k}x_{n+1-l},\endmultline$$
where $1\le j\le n$, and $m,n\in\Bbb N$, has a unique solution. Here
$\alpha_n^{ij}$ are the
coalgebra structure constants of $\Cal G_{\infty}$.
 The above system can be obtained
by differentiating (3b) with respect to $y$ and setting $y=e$, in
which case $\alpha^{ij}_n=\pd {\omega_{ij}}{y_n}\big|_{y=e}$. The
existence of a solution is furnished by Theorem 1 since any solution
of (3b) is a solution of the above system. To show
that it is unique one shows inductively that the corresponding
homogeneous system
$$\multline\sum_{i=j}^n(i+1-j)x_{i+1-j}\pd
{\omega_{mn}}{x_i}=\omega_{m+1-j,n}(x)(m+1-j)+\omega_{m,n+1-j}(x)(n+1-j)\endmultline$$
has only the trivial solution.
\par A subfamily of Lie bialgebra structures that corresponds to the family of
Poisson-Lie structures (11) is given by
$$\alpha_{d}(e_n)=2ne_d\wedge e_n-2(n-d)e_0\wedge e_{d+n},\tag 13$$
for each $d\in\Bbb N$. The entries of the $r$-matrix in this case are
$r^{ij}=\delta_{i+1}^{1}\delta_{j+1}^{d+1}-\delta_{i+1}^{d+1}\delta_{j+1}^{1}=\lambda_{i+1,j+1}$.
The family (13) of Lie bialgebra structures on
${\Cal G}_{\infty}$ had been found and studied in [2,7]. Also, we describe
below a 1-parameter family, $\alpha_{d,\lambda}$, of Lie bialgebra structures,
for each
$d\ge 1$, of which the family (13) is a subfamily obtained after the
specialization $\lambda=0$. Namely,
$$\eqalign{\alpha_{d,\lambda}(e_n)=&2\sum_{i=d+n}^{\infty}(2n-i){\lambda}^{i-(n+d)}(d-1)^{i-(n+d)}e_0\wedge e_i-2n\sum_{i=d}^{\infty}{\lambda}^{i-d}(d-1)^{i-d}e_i\wedge e_n\cr
&+{2}\sum_{i=d+n}^{\infty}\sum_{j=1}^{d-1}(2n-i)\lambda^{i+j-(n+d)}(d-1)^{i+j-(n+d+1)}e_i\wedge e_j\cr &+{2}\sum_{i=d}^{\infty}\sum_{j=n+1}^{d+n-1}(2n-j)\lambda^{i+j-(n+d)}(d-1)^{i+j-(n+d+1)}e_i\wedge e_j.}$$
Again, the right-hand-side of the above formula is an element of the completed
tensor
product $\Cal G_{\infty}\widehat{\otimes}\Cal G_{\infty}$. This family
corresponds to the following solution of (10):
$$\varphi_{d,\lambda }(u,v)=\frac {1}{[1-(d-1)\lambda
u][1-(d-1)\lambda v]}\Biggl\{ uv(v^d-u^d)+\lambda
du^2v^2(u^{d-1}-v^{d-1})\Biggr\}.$$
We conclude by noting that as a consequence of Theorem 4 the equation
(10) is a functional realization of the classical Yang-Baxter equation
for $\Cal G_{\infty}$.
\par The complete proofs of the above results will be published elsewhere
[5,6].

}


\Refs
\widestnumber\no{6}
\ref\no 1
\by V.G. Drinfel'd
\paper Hamiltonian Structures on Lie Groups, Lie Bialgebras and the Geometric
Meaning of the Classical Yang-Baxter Equations
\jour Soviet Math. Dokl. \vol 27 \yr 1983 \pages 68
\endref

\ref\no 2
\by W. Michaelis
\paper A Class of Infinite-Dimensional Lie Bialgebras Containing the Virasoro
Algebra
\jour Advances in Mathematics \vol 107 \yr 1994 \pages 365-392
\endref

\ref\no 3
\by M.A. Semenov-Tian-Shansky
\paper What is a Classical r-matrix?
\jour Funct. Anal. and its Applications
\vol 17 \yr 1983 \pages 259-272
\endref

\ref\no 4
\by M.A. Semenov-Tian-Shansky
\paper Dressing Transformations and Poisson Group Actions
\jour Publ. RIMS \vol 21 \yr 1985 \pages 1237-1260
\endref

\ref\no 5
\by O.S. Stoyanov
\paper Poisson-Lie Structures on Infinite-Dimensional Jet Groups and
Quantum Groups Related to Them \lang submitted
\endref

\ref\no 6
\by O.S. Stoyanov, Ph.D. Thesis
\publ Virginia Polytechnic Institute and State University, 1993
\endref

\ref\no 7
\by E.J. Taft
\paper Witt and Virasoro Algebras as Lie Bialgebras
\jour Journal of Pure and Applied Algebra, \vol 87
\yr 1993\pages 301-312
\endref
\endRefs
\bye